\documentclass[12pt,a4paper]{article}
\usepackage{amsmath,amssymb,amsfonts}
\def\Bbb{\mathbb}

\title{\bf Inversion of two cyclotomic matrices}

\author{Kurt Girstmair}
\date{}

\makeatletter
\let\@@maketitle=\maketitle
\def\maketitle{\def\thispagestyle##1{\relax}\@@maketitle}
\makeatother
\newtheorem{theorem}{Theorem}
\newtheorem{prop}{Proposition}
\newtheorem{lemma}{Lemma}

\def\BE{\begin{equation}}
\def\EE{\end{equation}}
\def\BD{\begin{displaymath}}
\def\ED{\end{displaymath}}
\def\BA{\begin{array}}
\def\EA{\end{array}}
\def\BEA{\begin{eqnarray*}}
\def\EEA{\end{eqnarray*}}
\def\BI{\bibitem}

\def\Z{\Bbb Z}

\def\XX{{\cal X}}
\def\RR{{\cal R}}

\def\phi{\varphi}
\def\EPS{\varepsilon}

\def\MB{\mbox}

\def\OV{\overline}

\def\WH{\widehat}

\def\DIV{\,|\,}
\def\NDIV{\, \nmid \,}

\def\MN{\medskip\noindent}
\def\STOP{\hfill$\Box$}

\begin{document}
\maketitle

\begin{abstract}

\noindent
Let $n\ge 3$ be a square-free natural number. We explicitly describe the inverses of the matrices
\BD
   (2\sin(2\pi jk^*/n))_{j,k} \enspace \MB{ and }\enspace (2\cos(2\pi jk^*/n))_{j,k},
\ED
where $k^*$ denotes a multiplicative inverse of $k$ mod $n$ and $j,k$ run through the set
$\{l; 1\le l\le n/2, (l,n)=1\}$. These results are based on the theory of Gauss sums.
\end{abstract}

\section*{1. Introduction and results}

In the paper \cite{Le} Lehmer states that there are only few classes of matrices for which explicit formulas for the determinant, the eigenvalues and the inverse are known. He gives
a number of examples of this kind. Further examples can be found in the papers \cite{Ca}, \cite{Mor} and \cite{Mol}. The closest analogue of the matrices considered here is contained in
the article \cite{Mol}, namely, the matrix
\BD
  (\sin(2\pi jk/n))_{j,k},
\ED
where $1\le j,k\le n, (jk,n)=1$. The author of \cite{Mol} determines the characteristic polynomial of this matrix, the multiplicities of the eigenvectors being quite involved. In the present note we describe the eigenvalues
of similar matrices $S$ and $C$.  The main results, however, are explicit formulas for the inverses $S^{-1}$ and $C^{-1}$ in the cases when these matrices are invertible. This is in contrast to the papers
we have quoted, since explicit formulas for inverses are scarcely given there.

Let $n\ge 3$ be a natural number. Let $\RR$ denote a system of representatives of the group $(\Z/\Z n)^{\times}/\{\pm 1\}$. Suppose that $\RR$ is ordered in some way. Typically,
$\RR$ is the set $\{k; 1\le k\le n/2, (k,n)=1\}$ with its natural order. For $k\in \Z$, $(k,n)=1$, let $k^*$ denote an inverse of $k$ mod n (so $kk^* \equiv 1$ mod $n$).
We define
\BD
   s_k=2\sin(2\pi k/n)\enspace \MB{ and }\enspace c_k= 2\cos(2\pi k/n),
\ED
where $k\in \Z$, $(k,n)=1$. We consider the matrices
\BD
  S=(s_{jk^*})_{j,k\in \RR} \enspace \MB{ and }\enspace C=(c_{jk^*})_{j,k\in \RR},
\ED
which we call the {\em sine matrix} and the {\em cosine matrix}, respectively.

We think that the matrices $S$ and $C$ deserve some interest not only  because of their simple structure but also by reason of their connection with cyclotomy, in particular,
with Gauss sums (see \cite{We} for the history of this topic).

In order to be able to enunciate our main results, we define
\BE
\label{1.2}
 \lambda(k)=|\{q; q\ge 3, q\DIV n, k\equiv 1 \MB{ mod } q\}|
\EE
for $k\in\Z$, $(k,n)=1$.
Furthermore, put
\BE
\label{1.2.1}
   \WH s_k=\frac 1n\sum_{l\in\RR} (\lambda(lk)-\lambda(-lk))s_l,
\EE
for $k\in\Z$, $(k,n)=1$.
For the same numbers $k$ put
\BE
\label{1.2.2}
   \WH c_k=\frac 1n\sum_{l\in\RR} (\lambda(lk)+\lambda(-lk)+\rho_n)c_l,
\EE
with
\BE
\label{1.3}
   \rho_n=\begin{cases} 2, & \MB{ if }  n \MB { is odd;}\\
                        4, & \MB{ if }  n \MB{ is even.}
          \end{cases}
\EE
Our main results are as follows.

\begin{theorem} 
\label{t1}
The sine matrix $S$ is invertible if, and only if, $n$ is square-free or $n=4$. In this case
\BD
  S^{-1}=(\WH s_{jk^*})_{j,k\in\RR},
\ED
with $\WH s_{jk^*}$ defined by {\rm (\ref{1.2.1})}.
\end{theorem} 

\begin{theorem} 
\label{t2}
The cosine matrix $C$ is invertible if, and only if, $n$ is square-free. In this case
\BD
  C^{-1}=(\WH c_{jk^*})_{j,k\in\RR}
\ED
with $\WH c_{jk^*}$ defined by {\rm (\ref{1.2.2})}.
\end{theorem} 

\MN
The entries of $S$ have the form $\pm s_l$, $l\in\RR$. This is due to the fact that
\BD
   s_{jk^*}=\EPS s_l
\ED
with $\EPS\in \{\pm 1\}$, $l\in\RR$, if $jk^*\equiv \EPS l$ mod $n$. In the same way we have
\BD
  \WH s_{jk^*}=\EPS \WH s_l
\ED
if $jk^*\equiv \EPS l$ mod $n$. This means that it suffices to compute the numbers $\WH s_l$ only for $l\in\RR$ in order to write down the matrix $S^{-1}$.
Indeed, this matrix arises from $S$ if we replace each entry $\EPS s_l$ of $S$ by the respective entry $\EPS \WH s_l$.

The same procedure works in the case of the cosine matrix, whose entries have the form $c_l$, $l\in\RR$.

\MN
{\em Example}. Let $n=15$ and $\RR=\{1,2,4,7\}$. Then $S$ can be written
\BE
\label{1.10}
 S= \left(\begin{array}{rrrr}       s_1& -s_7& s_4& -s_2\\
                                  s_2& s_1 & -s_7& -s_4\\
                                  s_4& s_2 &s_1   &s_7\\
                                   s_7& -s_4& -s_2 &s_1
 \end{array}\right).
\EE
Theorem \ref{t1} yields $\WH s_1=(3s_1-s_2+s_7)/15$, $\WH s_2=(-s_1-s_4-3s_7)/15$, $\WH s_4= (-s_2+3s_4+s_7)/15$, and $\WH s_7=(s_1-3s_2+s_4)/15$.
We obtain $S^{-1}$ if we put a circumflex on each $s$ occurring in (\ref{1.10}).

\MN
{\em Remark.}
If $n=p$ is a prime, Theorem \ref{t1} shows that $S^{-1}$ is particularly simple, namely, $S^{-1}=\frac 1p S^t$ ($S^t$ is the transpose of $S$). There is no analogue for the cosine matrix. For instance,
if $p=7$ and $\RR=\{1,2,3\}$, we have $\WH c_1=(3c_1+2c_2+2c_3)/7$. The prime number case of the sine matrix can also be settled by means of a simple trigonometric argument.
This, however, seems to be hardly possible if $n$ consists of at least two  prime factors $p>q\ge 3$.

\section*{2. Proofs}

First we prove Theorem \ref{t1}, then we indicate the changes required by the proof of Theorem \ref{t2}.
Let $\XX$ denote the set of Dirichlet characters mod $n$, and $\XX^-$ and $\XX^+$ the subsets of odd and even characters, respectively. The matrix $S$ is connected with $\XX^-$, whereas
$C$ is connected with $\XX^+$. We note the orthogonality relation
\BE
\label{2.2}
  \sum_{\chi\in\XX^-}\chi(k)=\begin{cases} 0,         & \MB{ if } k\not\equiv \pm 1 \MB{ mod }n;\\
                                           \phi(n)/2, & \MB{ if } k\equiv 1 \MB{ mod }n;\\
                                           -\phi(n)/2, & \MB{ if } k\equiv -1 \MB{ mod } n,
                              \end{cases}
\EE
see \cite[p. 210]{Ha}. Here $(k,n)=1$ and $\phi$ denotes Euler's function.

Suppose that the set $\XX^-$ is ordered in some way. Then we can define the matrix
\BD
X=\sqrt{n/\phi(n)}(\chi(k))_{k\in\RR,\chi\in\XX^-}.
\ED
Since $|\RR|=|\XX^-|=\phi(n)/2$, $X$ is a square matrix. We note the following lemma.

\begin{lemma} 
\label{l1}

The matrix $X$ is unitary, i.e., $X^{-1} =\OV X^t$ (the transpose of the complex-conjugate matrix).

\end{lemma} 

\MN
{\em Proof.} This is an immediate consequence of the orthogonality relation (\ref{2.2}) (observe that $\OV{\chi}(k)=\chi(k^*)$).
\STOP

Let $\zeta_n=e^{2\pi i/n}$ be the standard primitive $n$th root of unity.
For $\chi\in\XX^-$ let
\BE
\label{2.3}
  \tau(\chi)=\sum_{k=1}^n\chi(k)\zeta_n^k
\EE
the corresponding Gauss sum, see \cite[p. 445]{Ha}.  We consider the diagonal matrix
\BD
 T=\MB{diag}(\tau(\OV{\chi}))_{\chi\in\XX^-}.
\ED

\begin{prop} 
\label{p1}
The sine matrix $S$ is normal. Indeed,
\BD
 \OV X^tSX=-iT.
\ED

\end{prop} 

\MN
{\em Proof.} We show $XT\OV X^t=iS$. Obviously, the entry $(XT\OV X^t)_{j,k}$ equals
\BD
 \frac 2{\phi(n)}\sum_{\chi\in\XX^-}\chi(j)\tau(\OV{\chi})\OV{\chi}(k)=\frac 2{\phi(n)}\sum_{\chi\in\XX^-}\chi(jk^*)\sum_{(l,n)=1}\chi(l^*)\zeta_n^l,
\ED
where the index $l$ satisfies $1\le l\le n$, $(l,n)=1$.
This can be written
\BD
   \frac 2{\phi(n)}\sum_{(l,n)=1}\zeta_n^l\sum_{\chi\in\XX^-}\chi(jk^*l^*).
\ED
Now the orthogonality relation (\ref{2.2}), together with $\zeta_n^l-\zeta_n^{-l}=is_l$, shows that this is just $is_{jk^*}$.
\STOP

In order to study the vanishing of the eigenvalues of $S$, we use the reduction formula
\BE
\label{2.4}
 \tau(\OV{\chi})=\mu\left(\frac n{f_{\chi}}\right)\OV{\chi}_f\left(\frac n{f_{\chi}}\right)\tau(\OV{\chi}_f),
\EE
see \cite[p. 448]{Ha}.
Here $\mu$ means the M\"obius function, $f_{\chi}$ the conductor of the character $\chi$, $\chi_f$ the primitive character belonging to $\chi$ (which is a Dirichlet character mod $f_{\chi}$)
and $\tau(\OV{\chi}_f)$ the Gauss sum
\BD
   \sum_{k=1}^{f_{\chi}}\OV{\chi}_f(k)\zeta_{f_{\chi}}^k.
\ED
Since
\BE
\label{2.6}
  \tau(\chi_f)\tau(\OV{\chi}_f)=-f_{\chi}
\EE
(see \cite[p. 269]{Ha}), formula (\ref{2.4}) shows when the eigenvalue $-i\tau(\OV{\chi})$ vanishes. We obtain the following result.

\begin{prop} 
\label{p2}

The matrix $S$ is invertible if, and only if, $n$ is square-free or $n=4$.

\end{prop} 

\MN{\em Proof.}
If $n$ is square-free, then $n/f_{\chi}$ is square-free and $(f_{\chi},n/f_{\chi})=1$. By (\ref{2.4}) and (\ref{2.6}), all Gauss sums $\tau(\chi)$ are different from $0$.
If $n=4$ and $\chi\in\XX^-$, then $f_{\chi}=4$ and $n/f_{\chi}=1$.

Conversely, suppose that $n$ is not square-free and different from $4$.
Then one of the following three cases occurs. There is a prime $p\ge 3$ such that $p^2\DIV n$, or $4p\DIV n$, or $8\DIV n$.
In the first and the second case there is a character $\chi\in\XX^-$ with $f_{\chi}=p$. Accordingly, $\chi_f(n/f_{\chi})=0$ or $\mu(n/f_{\chi})=0$.
In the third case there is a character $\chi\in\XX^-$ with $f_{\chi}=4$. Therefore, $\chi_f(n/f_{\chi})=0$.
\STOP

\begin{lemma} 
\label{l2}

Let $n$ be square-free or equal to $4$. For $k\in \Z$, $(k,n)=1$, we have
\BD
   \sum_{\chi\in\XX^-}\frac{\chi(k)}{f_{\chi}}=\frac{\phi(n)}{2n}\:(\lambda(k)-\lambda(-k)),
\ED
the $\lambda$'s being defined by {\rm (\ref{1.2})}.

\end{lemma} 

\MN{\em Proof}.
Obviously,
\BD
  \sum_{\chi\in\XX^-}\frac{\chi(k)}{f_{\chi}}=\sum_{d\DIV n}\frac 1d\sum_{\chi\in\XX^-\atop f_{\chi}=d}\chi(k).
\ED
M\"obius inversion gives
\BD
   \sum_{\chi\in\XX^-\atop f_{\chi}=d}\chi(k)=\sum_{q\DIV d}\mu\left(\frac dq\right)\sum_{f_{\chi}\DIV q}\chi(k).
\ED
Here we note that the characters $\chi\in \XX^-$ with $f_{\chi}\DIV q$ are in one-to-one correspondence with the odd Dirichlet characters mod $q$.
Indeed, if $\chi\in\XX^-$, one defines the Dirichlet character $\chi_q$ mod $q$ in the following way. If $(j,q)=1$, there is an integer $l$ with $(l,n)=1$ such that $l\equiv j$ mod $q$.
Then $\chi_q(j)=\chi(l)$, see \cite[p. 217]{Ha}.
Accordingly,
\BD
 \sum_{f_{\chi}\DIV q} \chi(k)=\sum_{\chi_q} \chi_q(k).
\ED
From (\ref{2.2}) we obtain
\BE
\label{2.8}
  \sum_{\chi_q}\chi(k)=\begin{cases} 0,         & \MB{ if } q\le 2 \MB{ or } q\ge 3 \MB{ and } k\not\equiv \pm 1 \MB{ mod }q;\\
                                     \phi(q)/2, & \MB{ if } q\ge 3 \MB{ and } k\equiv 1 \MB{ mod }q;\\
                                     -\phi(q)/2, & \MB{ if } q\ge 3 \MB{ and } k\equiv -1 \MB{ mod }q
                              \end{cases}
\EE
(observe that there are no odd characters $\chi_q$ if $q\le 2$).
Therefore, we have
\BD
  \sum_{\chi\in\XX^-\atop f_{\chi}=d}\chi(k)=\sum_{q\DIV d, q\ge 3\atop k\equiv\pm 1\mod q}\pm\mu\left(\frac dq\right)\frac{\phi(q)}2,
\ED
where the $\pm$ sign in the summand corresponds to the respective sign in the summation index. If we write $d=q\cdot r$, we have
\BD
  \sum_{\chi\in\XX^-}\frac{\chi(k)}{f_{\chi}}=\sum_{q\DIV n, q\ge 3\atop k\equiv \pm 1\mod q}\pm\frac{\phi(q)}2\sum_{r\DIV \frac nq}\frac{\mu(r)}{qr}.
\ED
Since
\BD
\sum_{r\DIV \frac nq}\frac{\mu(r)}{r}=\prod_{p\DIV\frac nq}\left(1-\frac 1p\right)=\frac{\phi(n/q)}{n/q}
\ED
we obtain
\BE
\label{2.10}
  \sum_{\chi\in\XX^-}\frac{\chi(k)}{f_{\chi}}=\sum_{q\DIV n, q\ge 3\atop k\equiv\pm 1\mod q}\pm\frac{\phi(q)}{2q}\cdot \frac{\phi(n/q)}{n/q}.
\EE
However, $n$ is square-free or equal to $4$, and so $\phi(q)\phi(n/q)=\phi(n)$. This implies
\BD
  \sum_{\chi\in\XX^-}\frac{\chi(k)}{f_{\chi}}=\frac{\phi(n)}{2n}(\lambda(k)-\lambda(-k)).
\ED
\STOP

\MN{\em Proof of Theorem \ref{t1}.} By Proposition \ref{p1}, $S^{-1}=iXT^{-1}\OV X^t$, which means that the entry $(S^{-1})_{j,k}$, $j,k\in\RR$, of $S^{-1}$ is given by
\BD
  (S^{-1})_{j,k}=\frac{2i}{\phi(n)}\sum_{\chi \in \XX^-}\chi(j)\tau(\OV{\chi})^{-1}\OV{\chi}(k).
\ED
From (\ref{2.4}) and (\ref{2.6}) we obtain
\BD
  \tau(\OV{\chi})^{-1}=\frac{\mu(n/f_{\chi})\chi_f(n/f_{\chi})\tau(\chi_f)}{-f_{\chi}}=-\frac{\tau(\chi)}{f_{\chi}}.
\ED
Therefore,
\BD
  (S^{-1})_{j,k}=\frac{-2i}{\phi(n)}\sum_{\chi \in \XX^-}\frac{\chi(jk^*)}{f_{\chi}}\tau(\chi).
\ED
Now (\ref{2.3}) yields
\BD
  (S^{-1})_{j,k}=\frac{-2i}{\phi(n)}\sum_{(l,n)=1}\zeta_n^l\sum_{\chi \in \XX^-}\frac{\chi(ljk^*)}{f_{\chi}}.
\ED
By Lemma \ref{l2},
\BD
  (S^{-1})_{j,k}=\frac{-2i}{\phi(n)}\sum_{(l,n)=1}\zeta_n^l\frac{\phi(n)}{2n}(\lambda(ljk^*)-\lambda(-ljk^*)).
\ED
Altogether, we have
\BD
  (S^{-1})_{j,k}=\frac{-i}{n}\sum_{(l,n)=1}\zeta_n^l(\lambda(ljk^*)-\lambda(-ljk^*)).
\ED
On observing that $s_l=-i(\zeta_n^l-\zeta_n^{-l})$, we obtain Theorem \ref{t1}.
\STOP

\medskip
The setting of the {\em proof of Theorem 2} is slightly different. Indeed, the unitary matrix $X$ is defined by
\BD
  X=\sqrt{n/\phi(n)}(\chi(k))_{k\in\RR,\chi\in\XX^+}.
\ED
The cosine matrix $C$ is normal, and $\OV X^tCX=T$, with $T=\MB{diag}(\tau(\OV{\chi}))_{\chi \in\XX^+}$. In this case it is easy to see that $T$ (and, hence, $C$) is invertible
if, and only if, $n$ is square-free.
Instead of (\ref{2.6}) we have
\BD
  \tau(\chi_f)\tau(\OV{\chi}_f)=f_{\chi}.
\ED
The analogue of Lemma \ref{l2} reads
\BE
\label{2.12}
  \sum_{\chi\in\XX^+}\frac{\chi(k)}{f_{\chi}}=\frac{\phi(n)}{2n}\:(\lambda(k)+\lambda(-k)+\rho_n)
\EE
with $\rho_n$ as in (\ref{1.3}). This is due to the fact that the counterpart of formula (\ref{2.8}) takes the form
\BD
   \sum_{\chi_q}\chi(k)=\begin{cases} 1,         & \MB{ if } q\le 2;\\
                                      0,         & \MB{ if } q\ge 3 \MB{ and }  k\not\equiv \pm 1 \MB{ mod }q;\\
                                     \phi(q)/2, & \MB{ if }  q\ge 3 \MB{ and }  k\equiv \pm 1 \MB{ mod }q.
                              \end{cases}
\ED
Accordingly, formula (\ref{2.10}) has the equivalent
\BD
  \sum_{\chi\in\XX^+}\frac{\chi(k)}{f_{\chi}}=\sum_{q\DIV n, q\ge 3\atop k\equiv\pm 1\mod q}\frac{\phi(q)}{2q}\cdot \frac{\phi(n/q)}{n/q}+\sum_{d\DIV n,\atop 2\NDIV d}\frac{\mu(d)} d,
\ED
which gives (\ref{2.12}).
Up to these differences, the proof follows the pattern of the proof of Theorem \ref{t1}.


\vspace{0.5cm}
\noindent
Kurt Girstmair            \\
Institut f\"ur Mathematik \\
Universit\"at Innsbruck   \\
Technikerstr. 13/7        \\
A-6020 Innsbruck, Austria \\
Kurt.Girstmair@uibk.ac.at

\end{document}